\newcount\secno
\newcount\prmno
\newif\ifnotfound
\newif\iffound

\def\namedef#1{\expandafter\def\csname #1\endcsname}
\def\nameuse#1{\csname #1\endcsname}

\long\def\ifundefined#1#2#3{\expandafter\ifx\csname
  #1\endcsname\relax#2\else#3\fi}
\def\hwrite#1#2{{\let\the=0\edef\next{\write#1{#2}}\next}}

\toksdef\ta=0 \toksdef\tb=2
\long\def\leftappenditem#1\to#2{\ta={\\{#1}}\tb=
\expandafter{#2}%
                                \edef#2{\the\ta\the\tb}}
\long\def\rightappenditem#1\to#2{\ta={\\{#1}}\tb=
\expandafter{#2}%
                                \edef#2{\the\tb\the\ta}}

\def\lop#1\to#2{\expandafter\lopoff#1\lopoff#1#2}
\long\def\lopoff\\#1#2\lopoff#3#4{\def#4{#1}\def#3{#2}}

\def\ismember#1\of#2{\foundfalse{\let\given=#1%
    \def\\##1{\def\next{##1}%
    \ifx\next\given{\global\foundtrue}\fi}#2}}

\def\section#1{\vskip1truecm
               \global\def\currenvir{section}
               \global\advance\secno by1\global\prmno=0
               {\bf \number\secno. {#1}}}

\def\subsection{\global\def\currenvir{subsection}
                \global\advance\prmno by1
                \smallskip \ind{ (\number\secno.\number\prmno) }}
\def\subsec{\global\def\currenvir{subsection}
                \global\advance\prmno by1
                { (\number\secno.\number\prmno)\ }}

\def\proclaim#1{\global\advance\prmno by 1
                {\bf #1 \the\secno.\the\prmno$.-$ }}

\long\def\th#1 \enonce#2\endth{%
   \vglue4pt\proclaim{#1}{\it
#2}\global\def\currenvir{th}\smallskip}

\def\bib#1{\rm #1}
\long\def\thr#1\bib#2\enonce#3\endth{%
\medbreak{\global\advance\prmno by 1\bf#1\the\secno.\the\prmno\ 
\bib{#2}$\!.-$ } {\it
#3}\global\def\currenvir{th}\smallskip}
\def\rem#1{\global\advance\prmno by 1
{\it #1} \the\secno.\the\prmno$.-$ }


\def\isinlabellist#1\of#2{\notfoundtrue%
   {\def\given{#1}%
    \def\\##1{\def\next{##1}%
    \lop\next\to\za\lop\next\to\zb%
    \ifx\za\given{\zb\global\notfoundfalse}\fi}#2}%
    \ifnotfound{\immediate\write16%
                 {Warning - [Page \the\pageno] {#1} No reference found}}%
                \fi}%
\def\ref#1{\ifx\labellist\empty{\immediate\write16
                 {Warning - No references found at all.}}
               \else{\isinlabellist{#1}\of\labellist}\fi}

\def\newlabel#1#2{\rightappenditem{\\{#1}\\{#2}}\to\labellist}
\def\labellist{}

\def\label#1{%
  \def\given{th}%
  \ifx\given\currenvir%
    {\hwrite\lbl{\string\newlabel{#1}{\number\secno.\number
\prmno}}}\fi%
  \def\given{section}%
  \ifx\given\currenvir%
    {\hwrite\lbl{\string\newlabel{#1}{\number\secno}}}\fi%
  \def\given{subsection}%
  \ifx\given\currenvir%
    {\hwrite\lbl{\string\newlabel{#1}{\number\secno.\number
\prmno}}}\fi%
  \def\given{subsubsection}%
  \ifx\given\currenvir%
  {\hwrite\lbl{\string%
\newlabel{#1}{\number\secno.\number\subsecno.\number
\subsubsecno}}}\fi
  \ignorespaces}

\newwrite\lbl
\def\begin{
\newlabel{def}{1.1}
\newlabel{ext}{1.3}
\newlabel{defK3}{1.4}
\newlabel{irr}{2.6}
\newlabel{sen}{3.2}
\newlabel{set}{3.3}
\newlabel{sel}{3.4}
\newlabel{dual}{3.5}
\newlabel{mm}{5.1}}


\magnification 1250
\pretolerance=500 \tolerance=1000  \brokenpenalty=5000
\mathcode`A="7041 \mathcode`B="7042 \mathcode`C="7043
\mathcode`D="7044 \mathcode`E="7045 \mathcode`F="7046
\mathcode`G="7047 \mathcode`H="7048 \mathcode`I="7049
\mathcode`J="704A \mathcode`K="704B \mathcode`L="704C
\mathcode`M="704D \mathcode`N="704E \mathcode`O="704F
\mathcode`P="7050 \mathcode`Q="7051 \mathcode`R="7052
\mathcode`S="7053 \mathcode`T="7054 \mathcode`U="7055
\mathcode`V="7056 \mathcode`W="7057 \mathcode`X="7058
\mathcode`Y="7059 \mathcode`Z="705A
\def\spacedmath#1{\def\packedmath##1${\bgroup
\mathsurround =0pt##1\egroup$}
\mathsurround#1
\everymath={\packedmath}\everydisplay={\mathsurround=0pt}}
\def\nospacedmath{\mathsurround=0pt
\everymath={}\everydisplay={} } \spacedmath{2pt}
\def\qfl#1{\buildrel {#1}\over {\longrightarrow}}
\def\rfl#1{\,\nospacedmath\mathop{\hbox to
12truemm{\rightarrowfill}}\limits^{\scriptstyle#1}\,}
\def\phfl#1#2{\normalbaselines{\baselineskip=0pt
\lineskip=10truept\lineskiplimit=1truept}\nospacedmath\smash 
{\mathop{\hbox to 8truemm{\rightarrowfill}}
\limits^{\scriptstyle#1}_{\scriptstyle#2}}}
\def\hfl#1#2{\normalbaselines{\baselineskip=0truept
\lineskip=10truept\lineskiplimit=1truept}\nospacedmath
\smash{\mathop{\hbox to
12truemm{\rightarrowfill}}\limits^{\scriptstyle#1}_{\scriptstyle#2}}}
\def\diagramme#1{\def\normalbaselines{\baselineskip=0truept
\lineskip=10truept\lineskiplimit=1truept}   \matrix{#1}}
\def\vfl#1#2{\llap{$\scriptstyle#1$}\left\downarrow\vbox to
6truemm{}\right.\rlap{$\scriptstyle#2$}}
\def\mono{\lhook\joinrel\mathrel{\longrightarrow}}
\def\ptimono{\lhook\joinrel\mathrel{\rightarrow}}
\def\iso{\vbox{\hbox to .8cm{\hfill{$\scriptstyle\sim$}\hfill}
\nointerlineskip\hbox to .8cm{{\hfill$\longrightarrow $\hfill}}}}

\font\eightrm=cmr8         \font\eighti=cmmi8
\font\eightsy=cmsy8        \font\eightbf=cmbx8
\font\eighttt=cmtt8        \font\eightit=cmti8
\font\eightsl=cmsl8        \font\sixrm=cmr6
\font\sixi=cmmi6           \font\sixsy=cmsy6
\font\sixbf=cmbx6\catcode`\@=11
\def\eightpoint{%
\textfont0=\eightrm \scriptfont0=\sixrm
\scriptscriptfont0=\fiverm
\def\rm{\fam\z@\eightrm}%
  \textfont1=\eighti  \scriptfont1=\sixi  \scriptscriptfont1=\fivei
  \def\oldstyle{\fam\@ne\eighti}\let\old=\oldstyle
\textfont2=\eightsy \scriptfont2=\sixsy \scriptscriptfont2=\fivesy
  \textfont\itfam=\eightit
  \def\it{\fam\itfam\eightit}%
  \textfont\slfam=\eightsl
  \def\sl{\fam\slfam\eightsl}%
  \textfont\bffam=\eightbf \scriptfont\bffam=\sixbf
  \scriptscriptfont\bffam=\fivebf
  \def\bf{\fam\bffam\eightbf}%
  \textfont\ttfam=\eighttt
  \def\tt{\fam\ttfam\eighttt}%
  \abovedisplayskip=9pt plus 3pt minus 9pt
  \belowdisplayskip=\abovedisplayskip
  \abovedisplayshortskip=0pt plus 3pt
  \belowdisplayshortskip=3pt plus 3pt 
  \smallskipamount=2pt plus 1pt minus 1pt
  \medskipamount=4pt plus 2pt minus 1pt
  \bigskipamount=9pt plus 3pt minus 3pt
  \normalbaselineskip=9pt
  \setbox\strutbox=\hbox{\vrule height7pt depth2pt width0pt}%
  \normalbaselines\rm}\catcode`\@=12

\newcount\noteno
\noteno=0
\def\up#1{\raise 1ex\hbox{\sevenrm#1}}
\def\note#1{\global\advance\noteno by1
\footnote{\parindent0.4cm\up{\number\noteno}\
}{\vtop{\eightpoint\baselineskip12pt\hsize15.5truecm\noindent
#1}}\parindent 0cm}

\def\pc#1{\tenrm#1\sevenrm}
\def\tx{\kern-1.5pt -}
\def\cqfd{\kern 2truemm\unskip\penalty 500\vrule height 4pt 
depth 0pt width 4pt}
\def\ind{\par\hskip 0.8truecm\relax}
\def\indp{\par\hskip 0.4truecm\relax}
\def\moins{\mathrel{\hbox{\vrule height 3pt depth -2pt width
6pt}}}
\def\Ker{\mathop{\rm Ker}\nolimits}

\def\Pic{\mathop{\rm Pic}\nolimits}
\frenchspacing
\input xypic
\vsize = 25.6truecm
\hsize = 16truecm
\voffset = -.7truecm
\parindent=0cm
\baselineskip15pt

\begin
\centerline{\bf  Fano threefolds and K3 surfaces}
\smallskip 
\centerline{Arnaud {\pc BEAUVILLE}}\vskip1truecm
{\bf Introduction}\smallskip 
\ind A smooth anticanonical divisor in a Fano threefold is a K3
surface, endowed with a natural polarization (the restriction of
the anticanonical bundle). The question we address in this note is:
which K3 surfaces do we get in this way? The answer turns out to
be very easy, but it does not seem to be well-known, so the Fano
Conference might be a good opportunity to write it down.
\ind To explain the result, let us consider a component ${\cal
F}_g$ of the moduli stack\note{The frightened reader may
replace ``stack" by ``orbifold" or even ``space"; in the latter
case the word ``smooth" in the Theorem below has to be taken
with a grain of salt.} of pairs $(V,S)$, where $V$ is a Fano
threefold of genus $g$ and
$S$ a smooth surface in the linear system $|K_V^{-1} |$. Let ${\cal
K}_g$ be the moduli stack of polarized K3 surfaces of degree
$2g-2$. By associating to
 $(V,S)$ the surface $S$ we get a morphism of stacks
$$s_g :{\cal F}_g\longrightarrow {\cal K}_g\
.$$
\ind We cannot expect $s_g$ to be generically surjective, at
least if our Fano threefolds  have $b_2>1$: indeed for each
$(V,S)$ in ${\cal F}_g$ the restriction map  
$\Pic(V)\rightarrow \Pic(S)$ is injective by the weak Lefschetz
theorem, and this is a constraint on the K3 surface $S$. This map
is actually a lattice embedding when we  equip  $\Pic(V)$
 with the scalar product $(L,M)\mapsto (L\cdot M\cdot
K_V^{-1} )$; it maps the element $K_V^{-1}$ of $\Pic(V)$ to
the polarization of $S$.
\ind To
take this into account, we fix a  lattice $R$  with a
distinguished element $\rho$ of square $2g-2$, and we consider
the moduli stack 
${\cal F}_g^R$ parametrizing pairs $(V,S)$ with a lattice
isomorphism $R\iso \Pic(V)$ mapping $\rho$ to
$K_V^{-1}$. Let ${\cal K}_g^R$ be the algebraic stack
parametrizing   K3 surfaces $S$ 
together with an
 embedding  of $R$ as a primitive sublattice of $\Pic(S)$,
mapping
$\rho$ to an ample class. We have as before a
forgetful morphism
$s_g^R:{\cal F}_g^R
\rightarrow {\cal K}_g^R$.
\smallskip {\bf
Theorem}$.-$ {\it The morphism $s_g^R :{\cal
F}_g^R\rightarrow {\cal K}_g^R$ is smooth and generically
surjective; its  relative dimension at $(V,S)$ is
$b_3(V)/2$}. 
\ind As a corollary, a general K3 surface with given Picard 
lattice $R$ and polarization class $\rho\in R$ is an anticanonical
divisor in a Fano threefold if and only if $(R,\rho)\cong
(\Pic(V),K_V^{-1})$ for some Fano threefold $V$. 
\ind The proof of the Theorem is given in \S 3, after some
preliminaries on deformation theory (\S 1) and construction of the
moduli stacks (\S 2). We give  some comments in \S 4, and in
\S 5  we discuss the analogous question for curve sections of K3
surfaces. \smallskip \ind We will work for simplicity over ${\bf
C}$, though part of the results remain valid over an arbitrary
algebraically closed field.
 
\def\tg{T_X\langle Y\rangle}
\def\tgh{T_V\langle S\rangle}
\section{A reminder on deformation theory} 
\ind In this section we will quickly review two well-known
results on deformation theory  that are needed for the proof.
The experts are encouraged to skip this part.
\ind   Let
$X$ be a  smooth  variety, $Y$ a closed, smooth
subvariety of
$X$. We denote by
$\tg\i T_X$ the subsheaf of vector fields which
are tangent to $Y$, and by $r:\tg\rightarrow T_Y$  the
restriction map. 
\th{Proposition}
\enonce
The infinitesimal deformations of
$(X,Y)$ are controlled by the sheaf $\tg$ {\rm (that is,
obstructions lie in $H^2(X,\tg)$, 
 first order deformations  are parametrized by
$H^1$ and infinitesimal automorphisms by $H^0$).}  The map
which associates to  a first order deformation of $(X,Y)$ the
corresponding deformation of $Y$ is the induced map $H^1(r):
H^1(X,\tg)\rightarrow H^1(Y,T_Y)$. 
\endth
\label{def}
\ind This 
can be  extracted, for instance, from [R], but in such a simple
 situation it is more direct to apply  Grothendieck's
theory, as explained in [Gi], VII.1.2.  Let us sketch briefly how this
works. Put
$X_\varepsilon=X\otimes_{\bf C}{\bf C}[\varepsilon]$ and
$Y_\varepsilon=Y\otimes_{\bf C}{\bf C}[\varepsilon]$, with
$\varepsilon^2=0$; let ${\cal A}_{X,Y}$ (resp. ${\cal A}_Y$) be
the   sheaf of local automorphisms of
$ Y_\varepsilon\i X_\varepsilon$ (resp. $Y_\varepsilon$) which
induce the identity modulo $\varepsilon$.
  According to ({\it loc. cit.}),  since the
deformations of $Y\i X$ (resp. $Y$) are locally trivial, they
 are controlled by the sheaf ${\cal
A}_{X,Y}$ (resp.
${\cal A}_Y$) (technically, these deformations form a gerbe, and
the sheaf
${\cal A}$ is a band for this gerbe). So we just have to identify
these sheaves. For ${\cal A}_Y$ this is classical: a  section
of ${\cal A}_Y$ over an open subset $U$ of $Y$ is given by an
algebra automorphism
\vskip-12pt$${\cal O}_U[\varepsilon]\longrightarrow {\cal
O}_U[\varepsilon]$$ which must be of the form ${\rm
I+\varepsilon\,\delta }$, where $\delta $ is a derivation of ${\cal
O}_U$; this gives a group isomorphism ${\cal A}_Y\cong T_Y$.
Similarly a local automorphism of $(X,Y)$ is given by a diagram
\vskip-15pt$$\diagramme{{\cal
O}_X[\varepsilon]&\hfl{I+\varepsilon D}{} &{\cal
O}_X[\varepsilon]&\cr
\vfl{}{}&&\vfl{}{}&\cr
{\cal O}_Y[\varepsilon]&\hfl{I+\varepsilon \,\delta  }{}&{\cal
O}_Y[\varepsilon]&,
}$$\vskip-8pt where $D$ and $\delta $ are local derivations of
${\cal O}_X$ and ${\cal O}_Y$. The commutativity of the diagram
means that
$D$, viewed as a vector field, is tangent to $Y$, and induces the
vector field $\delta $ on $Y$. This gives an isomorphism ${\cal
A}_{X,Y}\cong \tg$; the forgetful map ${\cal A}_{X,Y}\rightarrow
{\cal A}_Y$ maps $(D,\delta)$ onto $ \delta $, thus coincides
with $r:\tg\rightarrow T_Y$.\cqfd

\subsection Let now $X$ be a smooth variety and
$R$  a free, finitely generated submodule of $\Pic(X)$; we 
consider the deformation problem for
$(X,R)$. Choosing a basis for $R$ this amounts
to deform $X$ together with line bundles $L_1,\ldots, L_p$.
As above the deformations of a pair
$(X,L)$ are controlled by the sheaf of local
automorphisms of $(X\otimes_{\bf C}{\bf C}[\varepsilon],
L\otimes_{\bf C}{\bf C}[\varepsilon])$ inducing the identity
modulo $\varepsilon$; this is readily identified with the
sheaf
${\cal D}^1(L)$ of first order differential operators of $L$, the
map $(X,L)\mapsto [X]$ corresponding to the symbol map
$\sigma :{\cal D}^1(L)\rightarrow T_X$ (this is of course
classical). Therefore deformations of
$(X,L_1,\ldots,L_p)$ are controlled by the sheaf ${\cal
D}^1(R):={\cal D}^1(L_1)\times _{T_X}\ldots\times _{T_X}{\cal
D}^1(L_p)$, which appears as an extension
$$0\rightarrow {\cal O}_X^p\longrightarrow {\cal
D}^1(R)\longrightarrow T_X\rightarrow 0 \ .$$The extension
class  lies in $H^1(\Omega^1_X)_{}^p$,  its $i$\tx th component
being the Atiyah class  $c_1(L_i)\in H^1(X,\Omega^1_X)$.  In a
more intrinsic way this can be written as an extension 
$$0\rightarrow R^*\otimes_{\bf Z}{\cal O}_X\longrightarrow
{\cal D}^1(R)\longrightarrow T_X\rightarrow 0
\eqno\subsec$$ whose class in
$H^1(X,\Omega^1_X)\otimes_{\bf Z}R^*$ is deduced from the
map $c_1:R\rightarrow H^1(X,\Omega^1_X)$.\label{ext}
\ind Assume now that $X$ is a K3 surface. We have $H^1(X,{\cal
O}_X)=H^2(X,T_X)=0$, and choosing a non-zero holomorphic
2-form on $X$ defines an isomorphism
$H^2(X,{\cal O}_X)\iso{\bf C}$.  The extension (\ref{ext})
gives rise to an exact sequence
$$0\rightarrow H^1(X,{\cal D}^1(R))\longrightarrow
H^1(X,T_X)\qfl{\partial} R^*\otimes _{\bf Z} {\bf
C}\longrightarrow H^2(X,{\cal D}^1(R))\rightarrow 0
$$
where $\partial $ is the cup-product with the extension class;
that is, for
$\xi\in H^1(X,T_X)$ and $L$ a line bundle in
$R$, we have
$\langle\partial (\xi),L\rangle=\xi\cup c_1(L)$. In other words,
using Serre duality, $\partial $ is the transpose of the natural
map $c_1:R\otimes_{\bf Z}{\bf C}\rightarrow
H^1(X,\Omega^1_X)$. Since $c_1$ is injective, $\partial $ is
surjective, hence $H^2(X,{\cal D}^1(R))=0$ and $H^1(X,{\cal
D}^1(R))=\Ker\partial $. Therefore:
\th Proposition
\enonce Let $X$ be a $K3$ surface and $R$ a subgroup of
$\Pic(X)$. The infinitesimal deformations of $(X,R)$ are
unobstructed. The first order deformations are parametrized by
the orthogonal of $c_1(R)\i H^1(X,\Omega^1_X)$ in
$H^1(X,T_X)$.\cqfd
\endth\label{defK3}
\section{The stacks ${\cal K}_g^R$
and ${\cal F}_g^R$}

\subsection Let $V$ be a smooth Fano threefold. Recall that the
{\it genus} $g$ of $V$ is defined by the formula $2g-2=(K_V^{-1}
)^3$. If $S$ is a smooth K3 surface in the linear system
$|K_V^{-1}|$, the induced polarization  $L:=K_V^{-1}{}_{|S}^{}$
satisfies
$L^2=2g-2$, so that the  curves of $|L|$ have genus $g$.
\ind As explained in the introduction, we will consider
 $\Pic(V)$ as a lattice with the product $(L,M)\mapsto
(L\cdot M\cdot K_V^{-1} )$. 
\subsection The definition of the moduli stack ${\cal F}$ of
pairs
$(V,S)$  is straightforward: we start from the moduli stack
${\cal T}$ of Fano threefolds. Let  $f:{\cal
V}\rightarrow {\cal  T}$ be the   universal family; the
projective bundle ${\bf P}((f_*K_{{\cal V/T}})^*)$ parametrizes
pairs
$(V,S)$ with $S\in |K_V^{-1}| $, and we take for ${\cal F}$ the
open substack defined by the condition that $S$ is smooth. We
add the subscript $g$ when we restrict to pairs $(V,S)$ of genus
$g$.
\subsection The definition of the moduli stacks ${\cal K}_g^R$
and ${\cal F}_g^R$ is slightly more involved.
Let
$f:X\rightarrow B$  be a smooth, projective morphism of
noetherian schemes. Following [G], we denote by
$\underline{\Pic}_{X/B}$ the sheaf on $B$ (for the faithfully flat
topology) associated to the presheaf
$(B'\rightarrow B)\mapsto\Pic(X\times_BB')$. According to {\it
loc. cit.}, this sheaf is representable by a group scheme over $B$,
for which we will use the same notation.  If
$f$ has relative dimension $2$, the intersection product defines a
bilinear form 
$\underline{\Pic}_{X/B}\times
\underline{\Pic}_{X/B}\rightarrow {\bf Z}_B$; the same holds in
(relative) dimension $3$ by taking the intersection product with
$K_{X/B}^{-1} $. 
\ind Let $R$ be a lattice, with a distinguished element $\rho$.
The moduli stacks ${\cal F}_g^R$ and ${\cal K}_g^R$ are defined
as follows. An object of
${\cal F}_g^R$ over a scheme $B$ is a pair $(V,S)$ over $B$,
where $V\rightarrow B$ is a family of Fano threefolds, of genus
$g$, and $S\i V$ a family of K3 surfaces over $B$, together with
a lattice isomorphism $R_B\iso \underline{\Pic}_{V/B}$
mapping $\rho$ onto the class of $K_V^{-1} $. Similarly, an object
of ${\cal K}_g^R$ over  $B$ is a family
 $S\rightarrow B$  of polarized K3 surfaces, of
genus $g$, together with a lattice embedding $R_B\mono
\underline{\Pic}_{S/B}$ mapping $\rho$ onto the 
polarization class. 
\ind That ${\cal K}_g^R$ and ${\cal F}_g^R$ are indeed
algebraic stacks follows from the result of Grothendieck quoted
above. Consider for instance the universal family ${\cal
S}\rightarrow {\cal K}_g$ of K3 surfaces with a genus $g$
polarization. Then
$\underline{\Pic}_{{\cal S}/{\cal K}_g}$ is representable by an
algebraic stack, which is a group scheme over ${\cal K}_g$. 
Choosing a basis
$(e_0,\ldots,e_p)$ of $R$ with $e_0=\rho$, we realize ${\cal
K}_g^R$ as an open and closed substack of
$(\underline{\Pic}_{{\cal S}/{\cal K}_g})^p$.
\ind Associating
to a pair $(V,S)$ over $B$ the family $S\rightarrow B$ with the
induced polarization and the composite map 
$R_B\iso \underline{\Pic}_{V/B}\mono 
\underline{\Pic}_{S/B}$ defines a morphism of stacks 
$s_g^R:{\cal F}_g^R \rightarrow  {\cal K}_g^R$.
\subsection Let us say a few words about the lattice $R$. In order
for our moduli stacks to be non-empty, $R$ must be a sublattice
of the Picard group of  a K3 surface, containing a
polarization;  also it must be isomorphic to the Picard lattice of a
Fano threefold. Thus:
\indp $\bullet$ $R$ is even, of   signature $(1,r-1)$;
\indp $\bullet$ $R$ has rank $r\le 10$; if $r\ge 6$, it is
isomorphic to the Picard lattice of $S_{11-r}\times {\bf
P}^1$, where $S_d$ is the Del Pezzo surface of degree $d$.
\ind (The latter property follows from Theorem 2 in [M-M]). 

\subsection Since $R$ has signature $(1,r-1)$, 
 the 
orthogonal of $\rho$  is negative definite, and therefore the
group of automorphisms of $R$ fixing $\rho$ is finite. 
It follows that the forgetful maps
${\cal F}_g^R\rightarrow {\cal F}_g$ and
${\cal K}_g^R\rightarrow {\cal K}_g$ are (representable and)
finite. The former map is actually is an \'etale covering, because
for any family
$V\rightarrow B$ of Fano threefolds the sheaf
$\underline{\Pic}_{V/B}$ becomes trivial on an \'etale
covering of $B$.
\ind  As for
the stack ${\cal K}_g^R$, we have
\th Proposition
\enonce The stack ${\cal K}_g^R$ is smooth, irreducible, of
dimension $20-r$.
\endth
\ind The smoothness and dimension of  ${\cal K}_g^R$ follow
from Proposition \ref{defK3}; its irreducibility is a consequence
of  the theory of the period mapping. Let us recall briefly how this
works, following the exposition in [D], 4.1. Let $L$ be
an even unimodular lattice of signature $(3,19)$ (all such
lattices are isomorphic). We choose an embedding of $R$ as a 
primitive sublattice of $L$ (such an embedding is unique up to
an automorphism of $L$ by Nikulin's results, see [D], thm. 1.4.8).
We consider  {\it marked} K3 {\it surfaces of
type} $R$, that is,  K3 surfaces $S$ with  a lattice
isomorphism
$\sigma :L\iso H^2(S,{\bf Z})$ such that $\sigma (R)$ is
contained in $\Pic(S)\i H^2(S,{\bf Z})$, and $\sigma
(\rho)$ is an ample class. These marked surfaces admit a fine
(analytic) moduli space $\widetilde{\cal K}_g^R$; the period map
induces an isomorphism of $\widetilde{\cal K}_g^R$ onto the
period domain $D_R$, which is the disjoint union of two
copies of a bounded symmetric domain of type IV ({\it loc. cit.}).
Our stack
${\cal K}_g^R$ is isomorphic to the quotient of 
$\widetilde{\cal K}_g^R$ by  the group $\Gamma _R$ of
automorphisms of $L$ which fix the elements of $R$. This group
acts on $D_R$  permuting its two connected components (this
can be seen exactly as in [B], Cor. p. 151). Thus the quotient stack
${\cal K}_g^R$ is irreducible.\cqfd\label{irr}

 \section{Proof of the theorem}
\subsection  By Proposition 
\ref{def} the infinitesimal behaviour of ${\cal F}_g$ (or ${\cal
F}_g^R$, since the forgetful map ${\cal F}_g^R\rightarrow
{\cal F}_g$ is \'etale) at a pair
$(V,S)$ is controlled by the sheaf $\tgh$, which is defined by the 
exact sequence
$$0 \rightarrow  \tgh\longrightarrow T_V
\longrightarrow N_{S/V} \rightarrow  0\
.\eqno\subsec$$\label{sen}  
We have 
$H^2(V,T_V)=H^2(V,\Omega^2_V\otimes K_V^{-1} )=0$ by the
Akizuki-Nakano theorem, and $H^1(S,N_{S/V})=0$ because
$N_{S/V}$ is an ample line bundle on $S$. Thus
 the exact sequence (\ref{sen})
gives $H^2(S, T_V\langle S\rangle)=0$, so that the first order
deformations of $(V,S)$ are unobstructed (in other words, the
stack ${\cal F}^R_g$ is smooth). 
\ind It follows from Proposition
\ref{def} that the tangent map to $s_g:{\cal F}_g\rightarrow
{\cal K}_g$ at
$(V,S)$  is
$H^1(r)$, where $r:\tgh\rightarrow T_S$ is the restriction map.
The map  $r$ is surjective, and its kernel is the subsheaf
$T_V(-S)$ of vector fields vanishing along
$S$, which in our case is isomorphic to $\Omega^2_V$. Thus 
we have an exact sequence
$$0\rightarrow  \Omega^2_V\longrightarrow \tgh
\qfl{r} T_S\rightarrow 0\ .\eqno\subsec$$
\ind Let us consider the associated  cohomology
exact sequence. 
 Since\label{set} 
$H^0(V,\Omega^2_V)$ and $H^0(S,T_S)$ are zero, we get first of
all 
$H^0(V,\tgh)=0$, so that $(V,S)$ has no
infinitesimal automorphisms (that is,   ${\cal F}^R_g$ is a
Deligne-Mumford stack). Then we get the
 exact sequence
$$0\rightarrow H^1(V, \Omega^2_V)\longrightarrow H^1(V,
T_V\langle S\rangle)\rfl{H^1(r)} H^1(S,T_S)\qfl{\partial} 
H^2(V,\Omega^2_V)\rightarrow 0\ .\eqno\subsec$$
\ind Let $i:S\ptimono V$ be the  inclusion map. To evaluate
$\partial$, consider the exact sequence\label{sel}
$$0\rightarrow  \Omega^1_V(\log S)(-S)\longrightarrow 
\Omega^1_V\qfl{i^*} \Omega^1_S\rightarrow
0\eqno\subsec$$\label{dual}
 deduced from (\ref{set}) by applying the duality functor
$R\underline{\rm Hom}_V(\ \ ,K_V)$ and using the canonical
isomorphisms $R\underline{\rm Hom}_V(T_S ,K_V)\cong 
R\underline{\rm Hom}_S(T_S ,K_S)\cong \Omega^1_S$.  By
general non-sense  the cohomogy exact sequence associated to
(\ref{dual}) is the dual of the one associated to (\ref{sel}); in
particular the transpose of $\partial$ is identified (through Serre
duality on
$V$ and  $S$) with the restriction map 
$H^1(i^*):H^1(V,
\Omega^1_V)\rightarrow H^1(S, \Omega^1_S)$ \hbox{-- up} to a
sign which is irrelevant for our purpose.
\ind Therefore $\Ker\partial$ is the orthogonal
of the image of $H^1(i^*)$. Because of the commutative diagram
$$\xymatrix{R\ar[d]\ar[dr]&\\
\Pic(V)\ar[r]^{i^*}\ar[d]^{c_1} & \Pic(S)\ar[d]\ar[d]^{c_1} \\
H^1(V,\Omega^1_V)\, \ar[r]^{H^1(i^*)}&\, 
H^1(S,\Omega^1_S)}$$
it is also the orthogonal of $c_1(R)\i
H^1(S,\Omega^1_S)$. By Proposition \ref{defK3} this is exactly
the tangent space to ${\cal K}_g^R$ at $S$, so the induced
map  $\tgh\rightarrow \Ker\partial $ is the tangent map to
$s_g^R$ at $(V,S)$. This proves that this  map  is surjective,
and the exact sequence (\ref{sel})
shows that its kernel is isomorphic to
$H^1(V,\Omega^2_V)$. Hence  $s_g^R$ is smooth, of relative
dimension $b_3(V)/2$, and generically surjective because ${\cal
K}_g^R$ is irreducible (Proposition \ref{irr}).\cqfd
\section{Consequences and comments} 
  \th Corollary 
\enonce Let $(S,h)$ be a polarized  $K3$ surface, $P$ its Picard
group; assume that $(S,h)$ is general  in ${\cal K}_g^P$.
Then $S$ is an anticanonical divisor in a Fano threefold if and only
if $(P,h)$ is isomorphic to
$(\Pic(V), K_V^{-1})$ for some Fano threefold $V$.\cqfd
\endth
\ind We leave to the reader
the enjoyable task of listing the pairs $(P,h)$ for the 87 types of
Fano threefolds with $b_2>1$ classified in [M-M]. In the case
$b_2=1$ we get the generic surjectivity of $s_g:{\cal
F}_g\rightarrow {\cal K}_g$; this is actually well-known, and
follows for instance from the work of Mukai [M1].
\subsection 
In most cases the map $s_g^R$ is not surjective.
 Consider for instance the component of
${\cal F}_5$ parametrizing  pairs $(V,S)$ with
$\Pic(V)={\bf Z}\cdot K_V$ and 
$g=5$. Each threefold $V$ is the complete intersection of 3
quadrics in
${\bf P}^6$, so we get in the image of $s_5$ all complete
intersections of 3 quadrics in
${\bf P}^5$,  which form a proper open
substack of ${\cal K}_5$ (it does not contain
hyperelliptic and trigonal K3 surfaces).

\subsection Part of the argument extends to Fano manifolds of
arbitrary dimension 
$n$, but the exact sequence (\ref{sel}) becomes
$$0\rightarrow H^1(V, \Omega^{n-1}_V)\longrightarrow H^1(V,
T_V\langle S\rangle)\longrightarrow H^1(S,T_S)\qfl{\partial} 
H^2(V,\Omega^{n-1}_V)\rightarrow 0\ ,$$so that the geometric
meaning of $\Ker \partial$ is not so clear. When $b_{n-1}(V)=0$
we see that the map $(V,S)\mapsto S$ is smooth.
\subsection A glance at the list of [M-M] shows that  roughly
half of the families of Fano threefolds have $b_3=0$; for these 
the map
$s_g^R$ is
\'etale, and one can ask whether it is an isomorphism 
onto an open substack. This is easy to prove in some cases
$(V={\bf P}^3, Q_3, {\bf P}^1\times {\bf P}^2,\ldots)$. For Fano
threefolds of index $2$ and genus $6$, it has been proved by
Mukai ([M1], Cor. 4.3). An interesting open case  is the one of
Fano threefolds of genus 12 with
$b_2=1$.
\section{K3 surfaces and canonical curves}
\subsection Let ${\cal KC}_g$ be the moduli stack of pairs $(S,C)$,
where $S$ is a  K3 surface  with a primitive polarization of
genus $g$, and
$C\i S$ a smooth curve  in the polarization class; let  ${\cal
M}_g$ be the moduli stack of curves of genus $g$. We have as
before a morphism of stacks
$$c_g:{\cal KC}_g\longrightarrow {\cal M}_g\ .$$ 
This morphism has been studied extensively. Let me
summarize  the main results. Recall first that
$\dim
{\cal KC}_g=19+g$ is greater than  
$\dim {\cal M}_g=3g-3$ for $g\le 10$, equal for $g=11$ and
smaller for
$g\ge 12$.
\ind $\bullet$ $c_g$ is generically surjective for $g\le 9$ and
$g=11$ [M1].
\ind $\bullet$ $c_g$ is  {\it not} surjective for $g=10$ [M1]; its
image is the hypersurface of ${\cal M}_g$ where the Wahl map
$\wedge^2H^0(C,K_C)\rightarrow H^0(C,K_C^{\otimes 3})$ fails
to be bijective [C-U].
\ind $\bullet$ $c_g$ is generically finite for $g=11$ and $g\ge
13$, but {\it not} for $g=12$ [M2].
\label{mm}
\subsection Let us consider the map $c_g$ from the differential
point of view that we have adopted in this note. Let $(S,C)\in{\cal
KC}_g$; we have by Serre duality $H^2(S, T_S\langle
C\rangle)=$ $H^0(S,\Omega^1_S(\log C))^*=0$, hence the
stack
${\cal KC}_g$ is smooth. By Proposition \ref{def}, the tangent
map to $c_g$ at $(S,C)$ is $H^1(r):H^1(S, T_S\langle C\rangle)
\rightarrow H^1(C,T_C)$. It appears in the  cohomology exact
sequence analogous to (\ref{sel})  
$$0\rightarrow H^1(S, T_S(-C))\longrightarrow H^1(S,
T_S\langle C\rangle)\ \hfl{H^1(r)}{}\ H^1(C,T_C)\qfl{\partial} 
H^2(S,T_S(-C))\rightarrow 0\ .$$
Using Serre duality, we see that $c_g$ {\it is smooth at $(C,S)$ if
and only if $H^0(S,\Omega^1_S(C))$ $=0$, and unramified at 
$(C,S)$ if and only if $H^1(S,\Omega^1_S(C))=0$}. Note that this
condition depends only on the polarization $L={\cal O}_S(C)$
and not on the particular curve $C$ in $|L|$ -- a fact which is not a
priori obvious.
\ind The results of (\ref{mm}) are thus equivalent to: 
\ind {\it Let $(S,L)$ be a general $K3$ surface with a 
 primitive polarization of genus $g$. We have: 
\ind $\bullet$ $H^0(S,\Omega^1_S\otimes L)=0$ for $g\le 9$
and $g=11$;
\ind $\bullet$ $ \dim H^0(S,\Omega^1_S\otimes L)=1$ for
$g=10$;
\ind $\bullet$ $H^1(S,\Omega^1_S\otimes L)=0$ for $g=11$
and} $g\ge 13$. 
\ind A direct proof of these results would provide an alternative
approach to the results of (\ref{mm}). 
\subsection Let us observe that though $c_g$ is generically
surjective for $g\le 9$ and
$g=11$, it  is {\it not} everywhere smooth. Take for instance a
K3 surface $S$ with an elliptic pencil $|E|$ and a smooth curve
$\Gamma $ of genus
$\gamma\in\{0,1\}$ with
$E\cdot\Gamma =2$; put
$L={\cal O}_S(kE+\Gamma) $. Then $L$ is a primitive polarization
of genus $2k+\gamma$. Let
$f:S\rightarrow {\bf P}^1$ be the map defined by the pencil $|E|$;
since $\Omega^1_S$ contains $f^*\Omega^1_{{\bf
P}^1}$, we get $\dim H^0(S,\Omega^1_S\otimes L)\ge k-1$. This
gives pairs $(S,C)$ in ${\cal KC}_g$, for $g\ge 4$, where
$c_g$ is not smooth. 
\ind Similarly, $c_g$  is {\it not} everywhere unramified for
$g=11$ or $g\ge 13$. A series of examples is provided by the
following result, which is essentially due to Mukai ([M2], Prop. 6):
\th Proposition
\enonce Let $V$ be a Fano threefold of index $1$ and genus $g$
such that
$K_V^{-1}$ is very ample, $S\in |K_V^{-1} |$  a 
$\!K3\!$ surface, $L:=K_V^{-1} {}_{|S}^{}$, $C$ a smooth
curve in the linear system $|L|$. The fibre of $c_g:{\cal
KC}_g\rightarrow {\cal M}_g$ at $(S,C)$ is positive-dimensional.
In particular, the space
$H^1(S,\Omega^1_S\otimes L)$ is non-zero.
\endth
{\it Proof} :  Consider $V$  embedded in ${\bf
P}(H^0(V,K_V^{-1} ))$. A general $C$ in $|L|$ is
contained in a Lefschetz pencil  $(S_t)_{t\in {\bf P}^1}$ of 
hyperplane sections of $V$: there is a finite subset $\Delta$ of
${\bf P}^1$ such that
$S_t$ is smooth for $t\in {\bf P}^1\moins \Delta$ and has an 
ordinary node for $t\in\Delta$.  The corresponding map $ {\bf
P}^1\moins\Delta\rightarrow {\cal K}_g$ goes to the boundary of
${\cal K}_g$ (consisting of K3 surfaces with a pseudo-polarization
of degree $2g-2$), and therefore cannot be constant. Thus we get a
1-dimensional family of pairs
$(S_t,C)$, for $t\in {\bf P}^1\moins \Delta$, which maps to the 
same point $[C]$ of ${\cal M}_g$.
 This  gives the result for $C$ general in
$|L|$, hence  for every smooth $C$ in $|L|$.\cqfd\smallskip 
\ind In view of the list in [M-M], we get examples of
positive-dimensional fibres of $c_g$ for all $g\le 28$ and
for $g=32$ (note that we want the polarization of $S$ to be
primitive, so
$V$ must be of index one). We  know no examples in  higher
genus, even with the weaker condition
$H^1(S,\Omega^1_S\otimes L)\not= 0$.

\vskip2cm
\centerline{ REFERENCES} \vglue15pt\baselineskip12.8pt
\def\num#1{\smallskip\item{\hbox to\parindent{\enskip [#1]
\hfill}}}
\parindent=1.3cm 
\num{B} A. {\pc BEAUVILLE}: {\sl Application aux espaces de
modules}.  G\'eom\'etrie des surfaces $K3$:
modules et p\'eriodes, exp. XIII. Ast\'erisque {\bf 126} (1985),
141--152.
\num{C-U} F. {\pc CUKIERMAN}, D. {\pc ULMER}: {\sl Curves of 
genus ten on $K3$ surfaces}. Compositio Math. {\bf 89} (1993),
 81--90.
\num{D} I. {\pc DOLGACHEV}: {\sl Integral quadratic forms:
applications to algebraic geometry (after V. Nikulin)}.
S\'em.   Bourbaki  1982/83, Exp. 611, 251--278. Ast\'erisque {\bf
105-106}, SMF, Paris (1983). 
\num{G} A. {\pc GROTHENDIECK}: {\sl Technique de descente et 
th\'eor\`emes d'existence en g\'eom\'etrie alg\'ebrique}. V.
   Les sch\'emas de Picard: th\'eor\`emes d'existence.  S\'em.
Bourbaki 1961/62, Exp.  232, 143--161.  SMF, Paris (1995). 
\num{Gi} J. {\pc GIRAUD}: {\sl Cohomologie non ab\'elienne}.
Grund.  math. Wiss. 
   {\bf 179}. Springer-Verlag, Berlin-New York (1971). 
\num{M1} S. {\pc MUKAI}: {\sl Curves, $K3$ surfaces and Fano
$\!3$\tx folds of genus}
$\leq 10$. Algebraic geometry and
   commutative algebra, Vol. I, 357--377, Kinokuniya, Tokyo
(1988).  
\num{M2} S. {\pc MUKAI}: {\sl  Fano $\!3$\tx folds}. Complex
projective geometry (Trieste--Bergen, 1989), 255--263,
   London Math. Soc. Lecture Note Ser. {\bf 179}, Cambridge Univ.
Press, Cambridge (1992). 
\num{M-M} S. {\pc MORI}, S. {\pc MUKAI}: {\sl Classification of
Fano $\!3$-folds with} $B\sb{2}\geq 2$. Manuscripta
   Math. {\bf 36} (1981/82), 147--162.
\num{R} Z. {\pc RAN}: {\sl Deformations of maps}. Algebraic
curves and projective geometry, 246--253. Springer Lecture Notes
{\bf 1389} (1989).   

\vskip1cm
\def\pc#1{\eightrm#1\sixrm}
\hfill\vtop{\eightrm\hbox to 5cm{\hfill Arnaud {\pc
BEAUVILLE}\hfill}
 \hbox to 5cm{\hfill Institut Universitaire de France\hfill}
\vskip-2pt\hbox to 5cm{\hfill \&\hfill}\vskip-2pt
 \hbox to 5cm{\hfill Laboratoire J.-A. Dieudonn\'e\hfill}
 \hbox to 5cm{\sixrm\hfill UMR 6621 du CNRS\hfill}
\hbox to 5cm{\hfill {\pc UNIVERSIT\'E DE}  {\pc NICE}\hfill}
\hbox to 5cm{\hfill  Parc Valrose\hfill}
\hbox to 5cm{\hfill F-06108 {\pc NICE} Cedex 02\hfill}}
\end